\documentclass[12pt]{article}
\usepackage[utf8]{inputenc}
\usepackage{amsmath}
\usepackage{amssymb}
\usepackage{amsthm}
\newcommand{\h}{\hspace{1cm}}
\newcommand{\z}{\varepsilon}
\newcommand{\ca}{ C^{\alpha, \frac{\alpha}{2}}(\bar{Q}_{T}) }
\newcounter{local}
\newcounter{locallocal}
\newcommand{\scl}{\stepcounter{local}}
\newcommand{\bys}{\begin{eqnarray*}}
\newcommand{\eys}{\end{eqnarray*}}
\newcommand{\by}{\begin{eqnarray}}
\newcommand{\ey}{\end{eqnarray}}
\begin{document}
\begin{center}

{\bf Asymptotic Analysis for a Nonlinear Reaction-Diffusion System \\ Modeling an Infectious Disease}

\bigskip
Hong-Ming Yin\footnote{Corresponding Author. Email: hyin@wsu.edu}\\
Department of Mathematics and Statistics\\
Washington State University\\
Pullman, WA 99164, USA.\\
and \\
Jun Zou\\
Department of Mathematics\\
The Chinese University of Hong Kong\\
Shatin, N.T., Hong Kong

\end{center}

\begin{abstract}
    In this paper we study a nonlinear reaction-diffusion system which models an infectious disease caused by bacteria such as those for cholera. One of the significant features in this model is that a certain portion of the recovered human hosts may lose  a lifetime immunity and could be infected again. Another important feature in the model is that the mobility for each species is allowed to be dependent upon both the location and time. 
    With the whole population assumed to be susceptible with the bacteria, the model is a strongly coupled nonlinear reaction-diffusion system. We prove that the nonlinear system has a unique solution globally in any space dimension under some natural conditions on the model parameters and the given data. Moreover, the long-time behavior and stability analysis for the solutions are carried out rigorously. In particular, we characterize the precise conditions on variable parameters about the stability or instability of all steady-state solutions. These new results provide the answers to several open questions raised in the literature. 

\end{abstract}


\newpage
\section{Introduction}

In biological, ecological, health and medical sciences,  researchers have a great deal of interest to establish  a suitable mathematical model for various infectious diseases. The current global pandemic attracts even more scientists to this field. There are many different mathematical models for an infectious disease in the literature. Roughly speaking, these models can be divided by two categories: a data-based discrete model and a continuous model based on a population growth (see \cite{GF2008,SR2013,WMI2018}). Our approach is based on a continuous model which provides a much more convenient tool to analyze the complicated dynamics of the interaction among susceptible, infected and recovered patients. A continuous model is typically governed by a system of ordinary differential equations (ODE model) or a system of partial differential equations (PDE model). For an ODE model, a monumental work was done in 1927 by Kermack and McKendrick \cite{KM1927}. Since then, a significant progress has been made in modeling and  analyzing various infectious diseases such as SIR, SEIR models and their various extensions. An ODE model often provides a clear and precise description of physical quantities and their relations. By using an ODE model, one can study detailed dynamical interaction between  viruses and various species as well as other qualitative properties such as reproduction numbers. This type of ODE models is widely adopted and used by researchers in all fields, particularly those in biological and health sciences. On the other hand, when one takes the movement of species across different geographical regions into consideration, it is necessary to include a diffusion process in a mathematical model to reflect the movement. This leads to modeling an infectious disease by using a system of partial differential equations (PDEs), often called reaction-diffusion equations. A well-known work \cite{HLBV1994} discussed a number of PDE models arising from biological, ecological and animal sciences and explained why the PDE approach is more appropriate in those areas. There are a large number of research studies, conference proceedings and monograph in both PDE and ODE models in the literature. We list only some of them here as examples, e.g., \cite{AM1979,MA1979,DG1999,DW2002,KP1994} for the SIR ODE models and   \cite{ABLN2007,FMWW2019,HLBV1994,JDH1995,LN1996,SLX2019} for the SIR PDE models.  Many more references can be found in a SIAM Review paper by Hethcote \cite{H2000} and the monograph by Busenberg and Cooke \cite{BC2012},  Cantres and Cosner \cite{CC1981}, Daley and Gani \cite{DG1999}, Lou and Ni \cite{LN1996}, etc. It is worth noting from the mathematical point of view that the PDE models present significant more challenges for scientists to study the dynamics of the solutions and to analyze qualitative properties of the solutions. Many important mathematical questions such as global existence and uniqueness are still open for some popular PDE models. 
This is one of the motivations for the current study.

In this paper we consider a mathematical model in a heterogeneous domain for an infectious disease caused by bacteria such as Cholera  without lifetime immunity. Without considering the diffusion-process of the population, the ODE models have been studied extensively  (see, e.g., \cite{AB2011,BC2012,DW2002,ESTD2002}). The model considered in this work is a direct extension of the ODE model. 
To describe the mathematical model, we introduce the following variables:
\begin{eqnarray*}
S(x,t) & = &  \mbox{Susceptible population concentration at location $x$ and time $t$}  \\
I(x,t) &  = &  \mbox{ Infected population concentration at location $x$ and time $t$}  \\
R(x,t) & = &  \mbox{ Recovered population concentration at location $x$ and time $t$}  \\
B(x,t) & = &  \mbox{Concentration of bacteria at location $x$ and $t$} 
\end{eqnarray*}.

We assume that the whole population is susceptible to the bacteria. Moreover, the rate of growth for the population, denoted by $b(x,t,S)$,  depends on location, time and the population itself.
 A classical example for $b$ is that the population growth follows a logistic growth model with a maximum capacity $k_1>0$:
\[ b(x,t,s)=b_0s(1-\frac{s}{k_1}),\]
where $b_0> 0$ represents the growth rate of the population.

The population reduction caused by infected patients is denoted by a nonlinear function $g_1(x,t,S,I,B)$ which is nonnegative.
A typical form of the nonlinear function $g_1$ is given by (see \cite{YW2017,Yin2020}):
 \[ g_1(x,t,S,I,B)=\beta_1SI+\beta_2Sh_1(B), ~~h_1(B)=\frac{B}{B+k_2},\]
 where $\beta_1, \beta_2$ are positive transmission parameters and $h_1(B)$ represents the maximum saturation rate of bacteria on
 human hosts and $k_2>0$.
 
The bacteria growth follows the same assumption, denoted by $g_2(x,t,s)$ with a maximum capacity $k_3>0$:
\[ g_2(x,t,s)=g_0s(1-\frac{s}{k_3}),\]
where $g_0> 0$ is the growth rate of the bacteria.

We also assume that the diffusion coefficients depend on location and time.
By extending the ODE model (see \cite{BC2012,DW2002,LW2011} etc.,), we obtain  the following reaction-diffusion system:
\setcounter{section}{1}
\setcounter{local}{1}
 \begin{eqnarray}
 S_t-\nabla\cdot [a_1(x,t)\nabla S] & = & b(x,t,S)-g_1(x,t,S,I,B)-d_1S+\sigma R,\\
 I_t-\nabla \cdot [a_2(x,t)\nabla I] & = & g_1(x,t,S,I,B)-(d_2+\gamma)I,\scl \\
 R_t-\nabla \cdot [a_3(x,t)\nabla R] & = & \gamma I-(d_3+\sigma)R,\scl \\
 B_t-\nabla \cdot [a_4(x,t)\nabla B] & = & \xi I +g_2(x,t,B)-d_4 B.\scl
 \end{eqnarray}
 
 The biological meaning of various parameters and functions in the model are given below (see \cite{ESTD2002,WW2015,YW2016}): 
 \begin{eqnarray*}
 a_i & = & \mbox{the diffusion coefficients, $i=1,2,3,4$},\\
 \gamma & = & \mbox{the recovery rate of infectious individuals},\\
 \sigma & = & \mbox{the rate at which recovered individuals lose immunity},\\
 d_i & = & \mbox{the natural death rate of species or bacteria},\\
 \xi & = & \mbox{the shedding rate of bacteria by infectious human hosts}.
 \end{eqnarray*}

 To complete the mathematical model, we assume that the system (1.1)-(1.4) holds in $Q_T=\Omega\times (0,T]$ for any $T>0$, where $\Omega$ is a bounded domain in $R^n$ with $C^2$-boundary $\partial \Omega$. The initial concentrations for all species are known and we assume that no species can cross the boundary $\partial \Omega$. This leads to the following initial and boundary conditions:
 \begin{eqnarray}
 & & (\nabla_{\nu}S,\nabla_{\nu}I,\nabla_{\nu}R, \nabla_{\nu}B) = 0, \h (x,t)\in \partial \Omega\times (0,T],\scl\\
 & & (S(x,0),I(x,0),R(x,0),B(x,0))=(S_0(x),I_0(x),R_0(x),B_0(x)),  x\in \Omega,\scl
 \end{eqnarray}
where $\nu$ represents the outward unit normal on $\partial \Omega$.

  We would like to give a short review about the known results for the above model. For the ODE system corresponding to (1.1)-(1.4), there are many studies for various interesting mathematical problems such as global existences, dynamical interaction between the bacteria and species (see, e.g.,  \cite{AB2011,DW2002,ESTD2002,TH1992}). The stability analysis is also carried out by several researchers (see  \cite{LW2011,SD2012,TW2011} etc.). When the movement of species is considered in the model, the corresponding PDE system is much more complicated to study. This is due to the fact that the maximum principle can not be applied for a system of reaction-diffusion equations. 
  It is a challenge to establish the global well-posedness for the PDE system (1.1)-(1.6). Nevertheless, when the space dimension is equal to $1$, under certain conditions on $g_1$ and $g_2$, the global existence is established (see \cite{WW2015,YW2016,YW2017,Y2018}). The reason is that the total population is bounded in $L^1(\Omega)$, which implies a global boundedness for $S(x,t)$ by using Sobolev embedding for the space dimension $n=1$. However, this method does not work when the space dimension $n$ is greater than $1$. In a SIAM review article
  (\cite{PS2000}), the authors considered the following system 
  (with $a$ and $b$ being two positive constants):
  \begin{eqnarray*}
  & & u_t-a\Delta u =f(u,v), \h x\in \Omega, ~t>0,\\
  & & v_t-b\Delta v =g(u,v), \h x\in \Omega, ~t>0,
  \end{eqnarray*}
  subject to appropriate initial and boundary conditions.
  Suppose $f(0,v), g(u,0)\geq 0$ for all $u,v \geq 0$. Then under the condition that 
  \[ f(u,v)+g(u,v)\leq 0,\]
  the $L^1$-norms of the nonnegative solutions $u$ and $v$ are bounded, i.e.,
  \[\sup_{t>0}\int_{\Omega} (u+v)dx\leq C.\]
   However, the solution $(u,v)$ may blow up in finite time when the space dimension is greater than 1 if no additional
   conditions on $f(u,v)$ and $g(u,v)$ are made. Therefore, as indicated in \cite{PS2000}, one must impose some additional conditions in order to obtain a global bound for a reaction-diffusion system. There are some interesting results for a general reaction-diffusion system when leading coefficients are constants. In 2000, under certain additional conditions, Pierre-Schmitt (\cite{PS2000}) introduced a dual method to establish such a bound for the reaction-diffusion system. In 2007, Desvillettes-Fellner-Pierre-Vovelle introduced in \cite{DFPV2007} an entropy condition originated by Kanel in 1990
   (\cite{KA1990}) and extended the dual method to a more general reaction-diffusion system with constant diffusion coefficients and established the global bound with a quadratic-growth reaction as long as a total mass is controlled ($L^1-$boundedness). In 2009, Caputo-Vasseur  \cite{CV2009} extended the entropy method to establish a global existence for a reaction-diffusion system where the nonlinear reaction terms grow at most sub-quadratically. One can see an interesting review by M. Pierre in 2010 \cite{PI2010}. 
   Caceres-Canizo extended in 2017 \cite{CC2017} to the case where the reaction terms grow at most quadratically under certain conditions on the steady-state solutions. In 2018, Souplet \cite{SOU2018} established the global well-posedness for a reaction-diffusion system with quadratic growth in the reaction. Very recently, some considerable progress was made for a reaction-diffusion system by Fellner-Morgan-Tang in 2019 \cite{FMT2019} and Morgan-Tang in 2020 \cite{MT2020}. They are able to derive a global bound for the solution of a reaction-diffusion as long as the diffusion coefficients are smooth and nonlinear reaction terms in the system satisfy a condition called an intermediate growth condition, which replaces the entropy condition. Their approach is based on a combination of the dual method and the entropy method. In 2021, Fitzgibbon-Morgan-Tang-Yin \cite{FMTY2021} studied a very general reaction-diffusion system with a controlled mass and  nonsmooth diffusion coefficients.
   They established the global well-posedness for the system with at most a polynomial growth for reactions. Moreover, several interesting examples as applications arising from biological, health sciences and chemical reactive-flow were studied in the paper. Those results made a substantial progress for a general reaction-diffusion system with a controlled mass. However, due to the nonlinearity in Eq. (1.1), these results do not cover the nonlinear system (1.1)-(1.4), particularly, we do not have growth conditions here on $g_1$ with respect to $(s_1, s_2, s_3)$ for the global existence (see Theorem 2.1 in section 2).

  The purpose of this paper has twofold. The first purpose is to establish the existence of a global solution to the generalized system (1.1)-(1.6) in any space dimension, without any restriction on parameters nor growth conditions with respect to $s_i$ for $g_1$. This extends a result obtained by the first author in his recent work \cite{Yin2020}.  Our method in this paper is based on some key ideas developed in \cite{Yin2020}. The special structure of the system (1.1)-(1.4) will also play a key role. We shall also use various techniques from the theories of elliptic and parabolic equations (see \cite{EVANS,LI1996,LSU}). To derive an a priori bound, we use a crucial result for a linear parabolic equation in the Campanato-John-Nirenberg-Morrey space from \cite{Yin1997}, which extends  the DiGoigi-Nash's estimate with weaker conditions for nonhomogeneous terms. The other purpose of the current work is to present the stability analysis of all steady-state solutions, which was not addressed in \cite{Yin2020}. In particular, for the following classical choices of the growth model \cite{CC1981}: 
  \by
  & &  {\ }\hskip-0.8truecm b(x,t,S)=b_0S\left(1-\frac{S}{k_1}\right), ~g_1( x,t,S,I,B)=\beta_1SI+\beta_2Sh_1(B), ~h_1(B)=\frac{B}{B+k_2}\scl\\
  &  & {\ }\hskip-0.8truecm g_2(B)=g_0B\left(1-\frac{B}{k_2}\right),\scl
  \ey
 we are able to precisely describe what conditions are needed for a steady-state solution to be stable or unstable. Roughly speaking, we shall demonstrate that under the  conditions: 
 \[ d_1>b_0, ~~d_2\geq 0, ~~d_3\geq 0, ~~d_4>g_0,\]
 the steady-state solution is stable. On the other hand, if either $d_1<b_0$ or $d_4<g_0$, then we can choose a set of suitable values for parameters  $\sigma, \gamma, \beta_1$ and $\beta_2$ such that the steady-state solution is unstable.
 This implies that our stability conditions are optimal. This stability analysis provides some important guidance to practitioners and scientists in biological, ecological and health sciences.  
 
 The paper is organized as follows. In Section 2 we first recall some function spaces which are frequently used in the subsequent analysis, and then state our main results. In Section 3, we prove the first part of the main results on global solvability of the system (1.1)-(1.6) (Theorem 2.1 and Corollary 2.1). In Section 4 we focus on a general stability analysis and obtain the sufficient conditions on parameters which ensure the stability of a steady-state solution.
   In Section 5, for a set of concrete functions $b(x,t,s), h_1(s)$ and $g_2(x,t,s)$ we give precisely conditions on the model parameters, under which a steady-state solution is stable or unstable. Finally, some concluding remarks are given in Section 6.
   
   Throughout the paper, we shall use $C$, with or without subscript, for a generic constant depending only on the given data in the model, including the upper bound of the terminal time $T$, and it may take a different value at each occurrence.

 \section{ Preliminaries and Statement of Main Results}
 
 For reader's convenience, we recall some standard function spaces which will be used frequently in the subsequent analysis.
 
 For $\alpha\in (0,1)$, we denote by $C^{\alpha}(\bar{\Omega})$ (or $C^{\alpha, \frac{\alpha}{2}}(\bar{Q}_{T})$) the  H\"older space in which every function is H\"older continuous with respect to $ x$ (or $(x,t))$ with exponent $\alpha$ in $\bar{\Omega}$ (or $(\alpha,\frac{\alpha}{2})$ in $\bar{Q}_{T}$). For $T=\infty$, we write $Q_T=\Omega\times (0, T)$ as $Q=\Omega\times (0,\infty).$
 
 For $p\geq 1$ and a Banach space $V$ with norm $||\cdot||_v$, we define
 \[ L^p(0,T; V)=\{ F(t): t\in [0,T]\rightarrow V; ~ ||F||_{L^p(0,T;V)}<\infty\},\]
 equipped with the norm
 \[ ||F||_{L^{p}(0,T;V)}=\left(\int_{0}^{T} ||F||_v^p dt\right)^{\frac{1}{p}}.\]
 When $V=L^{p}(\Omega)$, we simply write
 $L^p(Q_{T})=L^{p}(0,T;L^{p}(\Omega))$, 
 with its norm as $||\cdot||_p$.
 
Sobolev spaces $W^{k,p}(\Omega)$ and $W_{p}^{k,l}(Q_{T})$ are defined the same as in the classical references (see, e.g., \cite{EVANS}). Let $V_2(Q_T)=\{ u\in C([0,T];W_{2}^{1,0}(\Omega)): ||u||_{V_{2}}<\infty\} $ (see \cite{LSU})
 equipped with the norm
 \[ ||u||_{V_{2}}=\max_{0\leq t\leq T}||u||_{L^{2}(\Omega)}+\sum_{i=1}^{n}||u_{x_{i}}||_{L^2(Q_{T})}.\]
 
We will also use the Campanato-John-Nirenberg-Morry space $L^{2,\mu}(Q_T)$, which is defined as a subspace of $L^2(Q_{T})$ with its norm given by 
\[||u||_{L^{2, \mu}(Q_{T})}=||u||_{L^{2}(Q_{T})}+[u]_{2, \mu, Q_{T}}<\infty,\]
where
\[ [u]_{2, \mu, Q_{T}}=\sup_{\rho>0, z_0\in Q_{T}}\left(\rho^{-\mu}\int_{Q_{\rho}(z_0)}|u-u_Q|^2dxdt\right)^{\frac{1}{2}},\]
with $z_0=(x_0,t_0), Q_{\rho}(z_0)=B_{\rho}(x_0)\times (t_0-\rho^2, t_0]$, and $u_Q$ representing the average of $u$ over $Q_{\rho}(z_0)$ for any $Q_{\rho}(z_0)\subset  Q_{T}$;
see Troianiello $\cite{T1987} $ for its detailed definition and properties.
An important fact of the space is that
$L^{2, \mu+2}(Q_{T})$ is equivalent to $\ca$ with $\alpha=\frac{\mu-n}{2}$ if $n<\mu\leq n+2$ (Lemma 1.19 in \cite{T1987}). We shall write the norm of $L^{2, \mu}(Q_{T})$ as 
$||u||_{2, \mu}$. 

 We first state the basic assumptions for the diffusion coefficients and the known data involved in our model (1.1)-(1.4).
 All other model parameters are assumed to be positive constants throughout this paper.
 One can easily extend the well-posedness results to more general system when those parameters are functions of $(x,t)$ as long as the basic structure of the system is preserved.
 
 \ \\
 {\bf H(2.1).}  Assume that $a_i\in L^{\infty} (Q)$. There exist two positive constants $a_0$ and $A_0$ such that 
 \[ 0<a_0\leq a_i(x,t)\leq A_0, \h (x,t)\in Q_{T}, ~i=1,2,3,4.\]
 {\bf H(2.2).} Assume that all initial data $U_0(x):=(S_0(x), I_0(x), R_0(x), B_0(x))$ are nonnegative on $\Omega$. Moreover, $\nabla U_0(x)\in L^{2,\mu_0}(\bar{\Omega})^4$ with $\mu_0\in (n-2,n)$.
 
 \  \\
 {\bf H(2.3)}. (a) Let $b(x,t,s), d_i(x,t,s)$ and $g_2(x,t,s)$ be measurable in $Q\times R^+$ and locally Lipschitz continuous  with respect to $s$, and 
 $0\leq b(x,t,0), ~d_i(x,t,0)\in L^{\infty}(Q)$. Moreover, it holds 
 for some $M>0$ that 
 \[ d_i(x,t,s)\geq d_0\geq 0, ~~b_s(x,t,s)\leq b_0,  \h (x,t,s)\in Q\times[M,\infty).\]
 (b) Let $g_1(x,t,s_1,s_2,s_3)$ be measurable in $Q\times (R^+)^3$ and nonnegative, differentiable with respect to $s_1,s_2,s_3$, and 
 \bys
 & & g_1(x,t,0,s_2,s_3)\geq 0, \h  s_2, s_3\geq 0,\\
 & &  g_2(x,t,0)\geq 0, ~~g_{2s}(x,t,s)\leq g_0,  \h (x,t,s)\in Q\times R^+.
 \eys
where $k_1, k_2$ and $k_3$ represent the maximum capacity of the general population, the infected population and the bacteria, respectively.
 
 For convenience, we set $U(x,t)=(u_1,u_2,u_3,u_4)$ to be a vector-valued function defined in $Q_T$, with 
 \[ u_1(x,t)=S(x,t), ~u_2(x,t)=I(x,t), ~u_3(x,t)=R(x,t), ~u_4(x,t)=B(x,t), \, ~(x,t)\in Q_T. \]
 %We use $U(x,t)=(u_1,u_2,u_3,u_4)$ to be a vector-valued function %defined in $Q_T$. 
 The right-hand sides of the equations (1.1)-(1.4) are denoted by $f_1(x,t,U)$, $f_2(x,t,U)$, $f_3(x,t,U)$ and $f_4(x,t,U)$, respectively.
 With the new notation, the system (1.1)-(1.6) can be written as the following reaction-diffusion system:
 \setcounter{section}{2}
 \setcounter{local}{1}
 \begin{eqnarray}
 & & u_{1t}-\nabla\cdot [a_1(x,t)\nabla u_1]=f_1(x,t,U), \h (x,t)\in Q_{T},\\
 & & u_{2t}-\nabla\cdot [a_2(x,t)\nabla u_2]=f_2(x,t,U), \h (x,t)\in Q_{T},\scl \\
& &  u_{3t}-\nabla\cdot [a_3(x,t)\nabla u_3]=f_3(x,t,U), \h (x,t)\in Q_{T},\scl\\
 & & u_{4t}-\nabla\cdot [a_4(x,t)\nabla u_4]=f_4(x,t,U), \h (x,t)\in Q_T,\scl 
 \end{eqnarray}
 subject to the initial and boundary conditions:
 \begin{eqnarray}
 & & U(x,0)=U_0(x):=(S_0(x),I_0(x),R_0(x),B_0(x)), \h x\in \Omega, \scl \\
 & & \nabla_{\nu}U(x,t)=0,           \h (x,t)\in \partial \Omega\times (0,T].\scl
 \end{eqnarray}
 
  We define 
 \[ X=V_2(Q_{T})\bigcap L^{\infty}(Q_{T}).\]
 \ \\
 {\bf Definition 2.1.} We say $U(x,t)\in X^4$ is a weak solution to the problem (2.1)-(2.6) in $Q_{T}$ if it holds 
 for all functions $\phi_k\in X$ with $\phi_{kt}\in L^2(Q_{T}), 
 \phi_k(x,T)=0$ on $\Omega$ for $k=1,2,3,4$:
 \begin{eqnarray*}
 & & \int_{0}^{T}\int_{\Omega}\left[-u_k\cdot \phi_{kt}+a_k\nabla u_k\cdot \nabla\phi_k\right]dxdt\\
& & = \int_{\Omega} u_k(x,0)\phi _k(x,0)dx+\int_{0}^{T}\int_{\Omega}f_k(x,t,U)\phi_k(x,t)dxdt. \scl
 \end{eqnarray*}
 %for any test function $\phi_k\in X$ with $\phi_{kt}\in L^2(Q_{T}), \phi_k(x,T)=0$ on $\Omega$ %for all $k=1,2,3,4$.

 \ \\
 {\bf Theorem 2.1.} Under the assumptions H(2.1)-H(2.3), the problem (2.1)-(2.6) has a  weak solution in $X$ and the weak solution is nonnegative and bounded in $Q_T$ for any $T>0$. 
 Moreover, it holds that $u_i(x,t)\in C^{\alpha, \frac{\alpha}{2}}(\bar{Q}_{T})$ for $i=1,2,3,4$.
 
 Under some additional conditions on $b$ and $g_2$, we can deduce an uniform bound of the weak solution to the problem (2.1)-(2.6) in $Q$. 
 We state such a result for the special case
 which is needed in the subsequent asymptotic analysis.
 
 \ \\
 {\bf Corollary 2.1.} Under the conditions H(2.1)-(2.2), we further assume 
 \[ b_s(x,t,s)-d \geq \lambda_0>0, ~~g_{2s}(x,s)-d_4\geq \lambda_0>0,\h (x,t,s)\in Q\times [0,\infty),\]
 and 
 \[ \int_{0}^{\infty}\int_{\Omega}b_0(x,t)dxdt<\infty.\]
 Then the weak solution of the problem (2.1)-(2.6) is bounded globally in $Q$.

 \ \\
 {\bf Remark 2.1.} The weak solution obtained in Theorem 2.1 may grow to infinity as $t\to \infty$ if there is no additional conditions imposed on $b(x,t,S), g_2(x,t,s)$ and $d_1(x,t,s), d_4(x,t,s)$. On the other hand, if one assumes that
 $g_1$ and $g_2$ grow at most in a polynomial power with respect to $s_i$, then one can verify that the conditions in \cite{FMTY2021}
 hold. Consequently, a global bound in $Q$ can be deduced.

 The next theorem states our main stability results for the steady-state solutions to the problem (2.1)-(2.6).
 
 \ \\
 {\bf Theorem 2.2.} Under the condition H(4.1) (see Section 4), a steady-state solution 
 is asymptotically stable if 
 \[ d_1>B_0, ~~d_4>G_0,\]
 and the parameters $\beta_1, \beta_2, \gamma, \sigma$ are appropriately small, where $B_0$ 
 and $G_0$ are constants which depend on the steady-state solution.
 
 It turns out that the conditions in Theorem 2.2 are almost necessary in order to ensure the stability of each steady-state solution. In Section 5, we will see that when $b(x,t,s), g_1$ and $g_2(x,t,s)$ are of the form in (1.7)-(1.8),
  then we have a very precise set of conditions for the model parameters to ensure the local stability or instability for each steady-state solution. To avoid repetitions, we state this result in Section 5, 
  since there are many specific cases we have to consider.

 \section{Global Solvability and Proof of Theorem 2.1}
 
 In this section we first derive some a priori estimates for a weak solution to the system (2.1)-(2.6), then show the existence of a unique weak solution. Finally, we establish the global boundedness and the H\"older continuity.

\setcounter{section}{3}
\setcounter{local}{1}
\ \\
{\bf Lemma 3.1} Under the assumptions H(2.1)-(2.2), a weak solution of the system (2.1)-(2.6) is nonnegative.

This is a well-known result since each $f_i(x,t,u_1,u_2,u_3,u_4)$ 
is quasi-positive for $i=1,2,3,4$, and is also locally Lipschitz continuous with respect to each $u_k$ for $k=1,2,3,4$.
Interested readers may refer to \cite{BLS2002} for a detailed proof.

Next we apply the energy method to derive an a priori estimate in the space $V_2(Q_T)$.
 \ \\
{\bf Lemma 3.2} Under the assumptions H(2,1)-(2.3), there exists a constant $C_1$ such that
\[\sum_{k=1}^{4}||u_k||_{V_{2}(Q_T)}\leq C_1.\]
%where $C_1$ depends only on known data and an upper bound of $T$.\\
{\bf Proof}. We multiply Eq.(1.1) by $u_1$ and integrate over $\Omega$ to obtain
\begin{eqnarray*}
& & \frac{1}{2}\frac{d}{dt}\int_{\Omega}u_1^2dx+a_0\int_{\Omega}|\nabla u_1|^2dx+\int_{\Omega} g_1u_1dx +d_0\int_{\Omega}u_1^2dx
\nonumber\\
& & \leq \int_{\Omega}b(x,t,u_1) u_1 dx+\sigma\int_{\Omega}u_1u_3dx\nonumber \\
& & \leq C\int_{\Omega}[1+u_1^2] dx+ C\int_{\Omega}[u_1^2+u_3^2]dx,\scl
\end{eqnarray*}
where we have used the assumption H(2.3)(a) at the second estimate. 

We can perform a similar energy estimate for Eq.(1.3) to deduce 
\begin{eqnarray*}
& & \frac{1}{2}\frac{d}{dt}\int_{\Omega}u_3^2dx+a_0\int_{\Omega}|\nabla u_3|^2dx 
 \leq \gamma \int_{\Omega}u_2u_3dx \leq C\int_{\Omega}[u_2^2 + u_3^2]dx.\scl
\end{eqnarray*}

In order to derive an estimate for $u_2$, we make use of the special structure of the system (2.1)-(2.4). To do so, we define 
\[ v(x,t)=u_1(x,t)+u_2(x,t), \h (x,t)\in Q.\]
Then it is easy to see that $v(x,t)$ satisfies
\setcounter{section}{3}
\setcounter{local}{1}
\begin{eqnarray}
 v_t-\nabla\cdot [d_2\nabla v]&=&\nabla \cdot [(d_1-d_2)\nabla u_1]+f_1(x,t,U)+f_2(x,t,U), ~~(x,t)\in Q_{T},\\
 \nabla_{\nu}v(x,t)&=&0, \h (x,t)\in \partial \Omega\times (0,T],\scl\\
 v(x,0)&=&S_0(x)+I_0(x), \h x\in \Omega.\scl
\end{eqnarray} 

We multiply Eq.(3.1)  by  $ v$ and then integrate over $\Omega$ to obtain
\begin{eqnarray*}
& &   \frac{1}{2}\frac{d}{dt}\int_{\Omega}v^2dx +a_0\int_{\Omega}|\nabla v|^2dx\\
& & =-\int_{\Omega} [(d_1-d_2)\nabla u_1\cdot \nabla v]dx+\int_{\Omega}v[f_1(x,t,U)+f_2(x,t,U)]dx\\
& & :=J_1+J_2.
\end{eqnarray*}

A direct application of the Cauchy-Schwarz's inequality implies
\[ |J_1|\leq \z\int_{\Omega}|\nabla v|^2dx +C(\z)\int_{\Omega}|\nabla u_1|^2 dx.\]
On the other hand, using the fact that 
\[ f_1(x,t,U)+f_2(x,t,U)=b(x,t,u_1)-d_1u_1+\sigma u_3-(d_3+\gamma)u_2,\]
we readily derive that 
\begin{eqnarray*}
|J_2| & = & |\int_{\Omega}v[f_1(x,t,U)+f_2(x,t,U)]dx|\\
& \leq &  C\int_{\Omega}[v (1+u_1+ u_3)]dx
 \leq  C+C\int_{\Omega}[v^2+u_1^2+u_3^2]dx.
\end{eqnarray*}
%where $C$ depends only on known data.

Now choosing $\z=\frac{a_0}{2}$, we can readily derive from the above estimates that 
\begin{eqnarray*}
& & \frac{d}{dt}\int_{\Omega}v^2dx +a_0\int_{\Omega}|\nabla v|^2dx
\leq C+C\int_{\Omega}[v^2+u_1^2+u_3^2]dx.
\end{eqnarray*}
By combining the above energy estimates for $u_1, v$ and $u_3$, we can further deduce 
\begin{eqnarray*}
&& \frac{d}{dt}\int_{\Omega}[u_1^2+v^2+u_3^2]dx+\int_{\Omega}[|\nabla u_1|^2+|\nabla v|^2+
|\nabla u_3|^3]dx\\
& & \leq C\int_{\Omega}[u_1^2+v_2^2+u_3^2]dx,
\end{eqnarray*}
then a direct application of Gronwall's inequality implies
\begin{eqnarray*}
& & \sup_{0<t<T}\int_{\Omega}[u_1^2+v^2+u_3^2]dx+\int_{0}^{T}\int_{\Omega}[|\nabla u_1|^2+|\nabla v|^2+
|\nabla u_3|^3]dxdt\\
& & \leq C+C\int_{\Omega}[S_0^2+I_0^2+R_0^2]dx.
\end{eqnarray*}
%where $C$ depends only on known data.

Noting that $v=u_1+u_2$, we can write 
\begin{eqnarray*}
& & \int_{\Omega}|\nabla v|^2 dx= \int_{\Omega}[|\nabla u_1|^2+|\nabla u_2|^2]dx+2\int_{\Omega}
[(\nabla u_1)\cdot (\nabla u_2)]dx.
\end{eqnarray*}
But using the Cauchy-Schwarz's inequality, we can see 
\begin{eqnarray*}
 \int_{\Omega}
[(\nabla u_1)\cdot (\nabla u_2)]dx & \leq & \z\int_{\Omega}|\nabla u_2|^2dx+C(\z)\int_{\Omega}|\nabla u_1|^2dx\\
    & \leq & \z\int_{\Omega}|\nabla u_2|^2dx+C(\z)\int_{\Omega}[u_1^2+u_3^2]dx.
\end{eqnarray*}
Using the above estimates and choosing $\z$ to be sufficiently small, 
we can obtain
\begin{eqnarray*}
& & \int_{\Omega}[u_1^2+u_2^2+u_3^2]dx+
+\int\int_{Q_{T}}[|\nabla u_1|^2+|\nabla u_2|^2+|\nabla u_3|^2]dxdt\\
& &\leq C+  C\int_{\Omega}[S_0^2+I_0^2+R_0^2]dx.
\end{eqnarray*}
%where $C$ depends only on known data.

For $u_4$, we note that
\[ h_2(x,t,u_4)u_4\leq k_0(u_4^2+1).\]
Then we can readily derive from Eq. (2.4) that 
\[ \frac{d}{dt}\int_{\Omega}u_4^2dx +a_0\int_{\Omega}|\nabla u_4|^2dx\leq C\int_{\Omega}[u_2^2+u_4^2]dx.\]
Now an integration over $(0, T)$ implies 
\begin{eqnarray*}
& & \sup_{0<t<T}\int_{\Omega} u_4^2dx+\int\int_{Q_{T}}|\nabla u_4|^2dxdt
\leq C+C\int_{\Omega}B_0^2dx+C\int\int_{Q_{T}}u_2^2 dxdt\\
& & \leq C+C\int_{\Omega}[S_0^2+I_0^2+R_0^2+B_0^2]dx.
\end{eqnarray*}
%where $C$ depends only on known data and the upper bound of $T$.
This proof of Lemma 3.2 is now completed.
\hfill Q.E.D.

In order to derive more a priori estimats,  we need a crucial result about the Camapanto-John-Nirenberg-Morrey estimate for a general parabolic equation. For reader's convenience, we state the result in detail here (see Lemma 3.3 below).
Consider the parabolic equation:
\begin{eqnarray}
& & u_t-Lu =\sum_{i=1}^{n}f_i(x,t)_{x_i}+f(x,t), \h (x,t)\in Q_{T},\scl\\
& & u(x,t)=0 ~~\mbox{or} ~~u_{\nu}(x,t)=0, \h (x,t)\in \partial \Omega \times (0,T],\scl\\
& & u(x,0)=u_0(x), \h x\in \Omega.\scl
\end{eqnarray}
where 
$Lu:=(a_{ij}(x,t)u_{x_{i}})_{x_j}+b_i(x,t)u_{x_i}+c(x,t)$
is an elliptic operator. We assume there are positive constants $A_1,A_2$ and $A_3$ such that  $A=(a_{ij}(x,t)_{n\times n}$ is 
a positive definite matrix that satisfies
\[ A_0|\xi|^2 \leq a_{ij}\xi_i\xi_j \leq A_1|\xi|^2, \h \xi\in R^n,\]
and 
\[ \sum_{i=1}^{n}||b_i||_{L^{\infty}(Q_{T})}+||c||_{L^{\infty}(Q_{T})}\leq A_2<\infty.\]

\ \\
{\bf Lemma 3.3}. (\cite{Yin1997}) Let $u(x,t)$ be a weak solution of the parabolic equation (3.8)-(3.10).
Let $u_0\in C^{\alpha}(\bar{\Omega})$ with $u_0(x)=0$ on $\partial \Omega$, and $\nabla u_0\in L^{2,\mu_0}(\Omega)$
for some $\mu_0\in (n-2,n)$.
Then for any $\mu\in [0,n)$, there exists a constant $C$ such that
\[ ||\nabla u||_{L^{2, \mu}(Q_{T})}\leq C[||\nabla u_0||_{L^{2, (\mu-2)^+}(\Omega)}+||f||_{L^{2, (\mu-2)^+}(\Omega)}+\sum_{i=1}^{n}
||f_i||_{L^{2, \mu}(Q_{T})}].\]
Moreover, it holds that $u\in L^{2, \mu+2}(Q_{T})$ and 
\[ ||u||_{L^{2, 2+\mu}( Q_{T})}\leq C[||\nabla u_0||_{L^{2, (\mu-2)^+}(\Omega)}+||f||_{L^{2, (\mu-2)^+}(\Omega)}+\sum_{i=1}^{n}||f_i||_{L^{2, \mu}(Q_{T})}\]
for a constant $C$ that depends only on $A_0, A_1,A_2, n$ and $Q_{T}$.

\ \\
{\bf Lemma 3.4} Under the assumptions H(2.1)-(2.3), the weak solution of (2.1)-(2.4) satisfies
\[\sum_{k=1}^{4}||u_k||_{\ca}\leq C(T).\]
%where $C(T)$ depends only on the known data and the upper bound of $T$.\\
{\bf Proof}. Let $\mu\in (0, n)$ be arbitrary. 
By Lemma 3.3, we have
\begin{eqnarray}
& & ||\nabla u_3||_{L^{2, \mu}(Q_{T})}\leq C[||\nabla R_0||_{L^{2,(\mu-2)^+}(\Omega)}+||u_2||_{L^{2, (\mu-2)^+}(Q_{T})}+
||u_3||_{L^{2, (\mu-2)^+}(Q_{T})}].\scl
\end{eqnarray}
%where $C$ depends only on known data and the upper bound of $T$.

On the other hand, we note that $v(x,t)=u_1(x,t)+u_2(x,t)$ satisfies the system (3.1)-(3.3), 
so we can apply Lemma 3.3 again to obtain 
\begin{eqnarray}
& &  ||\nabla v||_{L^{2, \mu}(Q_{T})}\leq C[||\nabla v_0||_{L^{2,(\mu-2)^+}(\Omega)}+\sum_{i=1}^{3}||u_i||_{L^{2, (\mu-2)^+}(Q_{T})}].\scl
\end{eqnarray}

To derive the $L^{2,\mu}$-estimate for $u_1$, we note that
\[ u_{1t}-\nabla [a_1(x,t)\nabla u_1]\leq b(x,t,u_1)-d_1u_1+\sigma u_{3}=[b_s(x,t,\theta)-d_1]u_1+b(x,t,0)+\sigma u_3,\]
where $\theta$ is the mean-value between $0$ and $u_1$.
Using the facts that $b_s(x,t,s)$ and $b(x,t,0)$ are bounded, we can use the same calculations as in Lemmas 3.2 and 3.3 to obtain
\[ ||\nabla u_1||_{L^{2,\mu}(Q_{T})}\leq C[||\nabla S_0||_{L^{2,\mu}(\Omega)}+||u_3||_{L^{2,\mu}(Q_{T})}].\]

Now we can combine the $L^{2,\mu}(Q_T)$-estimates for $u_1, v$ and $u_3$ 
and note that $v=u_1+u_2$ to obtain 
for any $\mu\in [0, n)$ that 
\begin{eqnarray}
& & \sum_{i=1}^{3}||\nabla u_i||_{L^{2, \mu}(\Omega)} \leq C[||\nabla U_0||_{L^{2,(\mu-2)^+}(\Omega)}
+\sum_{i=1}^{3}||u_i||_{L^{2, (\mu-2)^+}(Q_{T})}]+C.\scl
\end{eqnarray}
Using the fact that $u_i\in V_2(Q_{T})$, we derive  
for any $\mu_1\in [0,2)$ that 
\begin{eqnarray}
 \sum_{i=1}^{3}||\nabla u_i||_{L^{2, \mu_1}(Q_{T})}\leq C[\sum_{i=1}^{3}||\nabla u_{i0}||_{L^{2}(\Omega)}+1]. \scl
\end{eqnarray}

Now we can apply the interpolation theory for the parabolic equation (2.3) 
(see Lemma 2.6 in \cite{Yin1997} ) to further deduce 
\[ ||u_3||_{L^{2,\mu_1+2}(\Omega)}\leq C[||u_2||_{L^{2}(\Omega)}+||u_3||_{L^{2}(\Omega)}
+||\nabla u_3||_{L^{2}(Q_{T})}]+C.\]

Next we go back to the system (2.1)-(2.3) and apply the same process for $\mu_2=\mu_1+2$ 
to obtain 
\begin{eqnarray}
 & & \sum_{i=1}^{3}||\nabla u_i||_{2, \mu_2, Q_{T}} \leq C[\sum_{i=1}^{3}||\nabla u_{i0}||_{L^{2,\mu_2}(\Omega)}+\sum_{i=1}^{3}||u_i||_{2, (\mu_2-2)^+,Q_{T}} +C].\scl
\end{eqnarray}
%where $C$ depends only on known data.
Then after a finite number of steps, we can deduce for any $\mu\in (0,n)$ that 
\begin{eqnarray}
& & \sum_{i=1}^{3}||u_i||_{L^{2, \mu+2}(\Omega)}\leq C[\sum_{i=1}^{3}||u_i||_{L^2(\Omega)}+||\nabla u_{i0}||_{L^{2, (\mu-2)^2}(\Omega)}\nonumber\\
& & \leq C[\sum_{i=1}^{3}||u_i||_{L^{2}(Q_{T})}+\sum_{i=1}^{3}||\nabla u_{i0}||_{L^{2, (\mu-2)^+}]}.\scl
\end{eqnarray}
%where $C$ depends only on known data.

Now we apply the interpolation theory again (see Lemma 2.6 in \cite{Yin1997}) to derive 
\begin{eqnarray*}
& & \sum_{i=1}^{3}||u_i||_{2, \mu_0+4, Q_{T}}
 \leq C[\sum_{i=1}^{3}||u_i||_{L^{2}(Q_{T})}+\sum_{i=1}^{3}||\nabla u_{i0}||_{L^{2, \mu_0}}.\scl
\end{eqnarray*}
%where $C$ depends only on known data.
But noting that $\mu_0\in (n-2, n)$, we can then obtain by Lemma 1.19 in \cite{T1987} that
\[ \sum_{i=1}^{3}||u_i||_{\ca}\leq C,\]
for $\alpha =\frac{\mu_0+2-n}{2}$.
The proof of Lemma 3.5 is now completed.
\hfill Q.E.D. 

\ \\
{\bf Proof of  Theorem 2.1}.  First of all, by using the energy method  we see that the weak solution of (2.1)-(2.6) must be unique since the solution is bounded and $f_k$ is locally Lipschitz continuous with respect to $u_i$ for all $k, i\in \{1,2,3,4\}$. With the a priori estimates in Lemmas 3.1-3.4,  there are several approaches, such as the truncation method and Galerkin finite element method, to prove the desired result
(see, e.g., \cite{BLS2002,FMTY2021, Yin2020}). Here we choose a different approach, 
the bootstrap argument (see \cite{YCW2017}), for the proof. 
Let $T\in (0,\infty)$ be any fixed number,
it is easy to show that the system (2.1)-(2.6) has a unique local weak solution in $X$ in $Q_{T_0}$ for some small $T_0>0$. 
Let
\[ T^*=sup\{T_0: \mbox{the system (2.1)-(2.6) has a unique weak solution in $Q_{T_0}$}\}.\]
Suppose $T^*<T$ (otherwise, nothing is needed to prove). We note that the a priori estimates  in Lemmas 3.1 and 3.4 hold
for any weak solution. It follows that
\[ \lim_{t\rightarrow T*-}sup  [\sum_{k=1}^{4}||u_k||_{V_2(Q_t)}+\sum_{k=1}^{4}||u_k||_{\ca}]<\infty.\]
By the compactness, we know that
\[ u_k(x,T^*)\in H^1(\Omega), \nabla u_k\in L^{2, (\mu-2)^+}(\Omega) ~~\mbox{for any $\mu\in (n,n+2).$}\]

Now, we use $U(x,T^*)$ as an initial value and consider the system (2.1)-(2.6) for $t\geq T*$. Then the local existence result
implies that there exists a small $t_0>0$ such that the problem (2.1)-(2.6) has a unique weak solution in the interval $[T*, T^*+t_0).$
Consequently, we obtain a weak solution to the system (2.1)-(2.6) in the interval $[0,T^*+t_0)$.
This is a contradiction with the definition of $T^*$, therefore we have $T^*=T$.
\hfill Q.E.D.

\ \\

Next, we prove Corollary 2.1.
Assume that there exists a constant $\lambda_0>0$ such that
\[ d_1(x,t,s)-b_s(x,t,s)\geq \lambda_0>0, ~~d_4(x,t,s)-g_{2s}(x,t,s)\geq \lambda_0,\h (x,t,s)\in Q\times [0,\infty).\]
With the above assumption, we take the integration over $\Omega$ for Eq. (2.1)-(2.3) to obtain
\[ \frac{d}{dt}\int_{\Omega}(u_1+u_2+u_3)dx+\min\{(d_0,\lambda_0\}\int_{\Omega}(u_1+u_2+u_3)dx\leq \int_{\Omega}b(x,t,0)dx.\]
Then it is easy to see 
\[ \sup_{t\geq 0}\int_{\Omega}(u_1+u_2+u_3)dx\leq C.\]
%where $C$ depends only on known data.

Now we derive a uniform estimate in $L^2(Q)$. By using the energy estimate for Eq.(2.1), 
we can see that 
\[ \frac{d}{dt}\int_{\Omega}u_1^2dx+\int_{\Omega}|\nabla u_1|^2dx
\leq C[\int_{\Omega}b(x,t,0)^2dx+C\int_{\Omega}u_3^2dx.]\]
For $v(x,t):=u_1(x,t)+u_2(x,t)$, we can derive from Eq.(3.1)-(3.3) that
\bys
& & \frac{d}{dt}\int_{\Omega}v^2dx+\int_{\Omega}|\nabla v|^2dx
 \leq C\int_{\Omega}|\nabla u_1|^2dx+C\int_{\Omega}[(b(x,t,0)^2+u_1^2+u_3^2]dx\\
& & \leq C[\int_{\Omega}(b(x,t,0)^2+u_1^2+u_3^2) dx].
\eys
where we have used the estimate of $u_1$ at the second estimate.

Again, we can use the energy estimate for Eq.(2.3) to obtain
\[ \frac{d}{dt}\int_{\Omega}u_3^2dx+\int_{\Omega}|\nabla u_3|^2dx\leq C\int_{\Omega}u_2^2dx.\]
But we know from the Gagliardo-Nirenberg estimate 
for $p=q=2, s=1, \theta=\frac{n}{n-2}$ and $\varepsilon>0$, 
\[ \int_{\Omega} u^2dx\leq \varepsilon \int_{\Omega}|\nabla u|^2 dx+C(\varepsilon) ||u||_{L^1(\Omega)}, \]
%where $\varepsilon>0$.
then using the uniformly boundedness of $L^1(\Omega)$-norms of $u_1, u_2, u_3$, 
we get for sufficiently small $\varepsilon$, 
\bys
& & \sup_{t\geq 0}\int_{\Omega}[u_1^2+u_2^2+u_3^2]dx+\int_{0}^{t}\int_{\Omega}[|\nabla u_1|^2+|\nabla u_2|^2+|\nabla u_3|^2dx]\\
& & \leq C_1+C_2\int_{0}^{t}\int_{\Omega}b(x,t,0)^2dxdt
\leq C_3\,.
\eys
%where $C_3$ depends only on known data.

Next we use the iteration method again as in the proof of Theorem 2.1.  
From Eq.(3.2) for $v$ and $u_3$, we deduce, respectively, 
\[ ||\nabla v||_{2, \mu}\leq C+C||u_1||_{2, \mu}+C||u_3||_{2,\mu}\]
and
\[||\nabla u_3||_{2,\mu}\leq C+C||u_2||_{2, \mu}. \]
%where $C$ depends only on known data.

For $u_1$, we see by noting that $g_1\geq 0$, 
\[ u_{1t}-\nabla[a_1(x,t)\nabla u_1]\leq [b_0(x,t)-d_1]u_1+\sigma u_3.\]
As $u_1\geq 0$ in $Q$, we can follow the same argument as in \cite{Yin1997} to obtain
for $\mu\in (n-2, n)$, 
\[ ||\nabla u_1||_{2,\mu}\leq C+C[||u_1||_{2, \mu}+||u_3||_{2, \mu}].\]
%where $\mu\in (n-2, n)$ and $C$ depends only on the known data.

As $u_1, v, u_3$ are uniformly bounded in $L^2(Q)$, the interpolation for $v$ and $u_3$ with $\mu=0$ yields that
\[ ||v||_{2,2}+||u_3||_{2, 2}\leq C.\]
%where $C$ depends only on known data.
Hence, we can obtain the $L^{2, \mu}(Q)$-estimate for $\nabla u_1$ with $\mu=2$:
\[ ||\nabla u_1||_{2,2}\leq C+ C[||u_1||_{2, 2}+||u_2||_{2, 2}],\]
which is uniformly bounded. 

We can now go back to the equations for $v$ and $u_3$ with $\mu=2$ to obtain
\[  ||v||_{2,4}+||u_3||_{2, 4}\leq C[||u_1||_{2,2}+||u_2||_{2,2}+||u_3||_{2,2}].\]
%where $C$ depends only on known data.

By continuing the above iteration process, after a finite number of steps, 
we obtain for $\alpha=\frac{\mu_0-n}{2}$ that   
\[ ||v||_{\ca}+||u_3||_{\ca}\leq C.\]
%where $\alpha=\frac{\mu_0-n}{2}$ and $C$ depends only on know constants.
Consequently, we get
\[ ||u_1||_{L^{\infty}(Q)}\leq C.\]
%where $C$ depends only on known data.

Once we know that $u_2$ is uniformly bounded, then from Eq.(2.4), we 
can apply the maximum principle to obtain 
\[\sup_{t\geq 0}||u_4||_{L^{\infty}(\Omega)}\leq C.\]
%where $C$ depends only on known data.

With the a priori bound for each $u_i$, we can extend the weak solution in $Q_{T}$ to $Q$.

\hfill Q.E.D.

\section{Linear Stability Analysis}

To illustrate the main idea, we assume that $b$ and $g_2$ depend only on $x$ and $s$.
We also focus on the following model cases: 
\[ b(x,t,s)=b_0(x)s(1-\frac{s}{k_1}), ~~g_1=\beta_1u_1u_2+\beta_2 \frac{u_1u_4}{u_4+k_2}, ~~g_2=g_0(x)s(1-\frac{s}{k_3}).\]
Moreover, we assume that all parameters $\sigma, \gamma, \beta_1, \beta_2, d_i, k_i$ are positive constants. The general case can be carried out similarly as long as the functions are differentiable.

Consider the steady-state problem in $\Omega$:
\setcounter{section}{4}
\setcounter{local}{1}
\begin{eqnarray}
 -\nabla\cdot [a_1(x)\nabla u_1] & = & b(x,u_1)-g_1(x,u_1,u_2,u_4)-d_1u_1+\sigma u_3,\\
 -\nabla\cdot [a_2(x)\nabla u_2] & = & g_1(x,u_1,u_2,u_4)-(d_2+\gamma)u_2,\scl \\
 -\nabla\cdot [a_3(x)\nabla u_3] & = & \gamma u_2-(d_3+\sigma)u_3,\scl \\
 -\nabla\cdot [a_4(x)\nabla u_4] & = & \xi u_2 +g_2(x,u_4)-d_4 u_4\scl
 \end{eqnarray}
 subject to the boundary condition
 \begin{eqnarray}
 \partial_{\nu}U(x) =0, \h x\in \partial \Omega,\scl
 \end{eqnarray}
 where  $U(x)=(u_1(x), u_2(x), u_3(x), u_4(x))$.
 
 It is clear that there is a trivial solution $U(x)=(0,0,0,0)$ if $b(x,0)=g_1(x,0)=g_2(x,0)=0.$ But we are interested in nontrivial solutions, and will make the following assumptions.
 
 \ \\
 H(4.1). (a) $0<a_0\leq a_i(x)\leq A_0$ on $\Omega$; \\
         (b) $b_0(x)\geq b_1>0$ and $g_0(x)\geq g_1>0$, and both are bounded.
 
 \ \\
 {\bf Lemma 4.1}. Under the assumptions H(4.1), the elliptic system (4.1)-(4.5) has at least one nonnegative weak solution
 $U(x)\in W^{1,2}(\Omega)$. Moreover, the weak solution is H\"older continuous in $\bar{\Omega}$ for any space dimension.\\
 {\bf Proof.}  Since the argument is very similar to the case for a parabolic system, we only sketch the proof.
 The key step is to derive an a priori estimate in H\"older space. 
 As a first step, we know that a solution of (4.1)-(4.5) must be nonnegative since every right-hand side of (4.1) to (4.4) is quasi-positive. Next we can use the same argument as for the parabolic case to
 derive $L^1$-estimate for $u_i(x)\geq 0, i=1,2,3,4$ on $\Omega$.
 Indeed,  by direct integration we have 
 \[ \int_{\Omega}[d_2u_2+d_3u_3]dx+\int_{\Omega}b_0(x) u_1^2dx=\int_{\Omega}(b_0-d_1)u_1dx.\]
 Then an application of the Cauchy-Schwarz's inequality yields
 \[ \int_{\Omega}[u_1^2+d_2u_2+d_3u_3]dx\leq C.\]
 %where $C$ depends only on $b_0, d_1$ and $\Omega$.
 
 On the other hand, we obtain from Eq.(4.1) that 
 \[ g_0\int_{\Omega}u_4^2dx \leq \xi\int_{\Omega}u_2dx+g_0\int_{\Omega}(g_0-d_4)u_4dx\leq C+\frac{g_0}{2}\int_{\Omega}u_4^2dx,\]
 which implies 
 \[ \int_{\Omega}u_4^2dx \leq C.\]
 %where $C$ depends only on known data.
 
 Next step is to derive the $L^2(\Omega)$-estimate for $u_2$ and $u_3$. The idea is very much similar to the case for a parabolic system.
 The energy estimate for Eq.(4.1) yields that, for any $\varepsilon>0$,
 \[ \int_{\Omega}|\nabla u_1|^2dx+\int_{\Omega}u_{1}^3dx\leq C(\varepsilon)+\varepsilon\int_{\Omega}u_3^2dx.\]
 It is easy to see that, by adding up Eq.(4.1) and Eq.(4.2), $v(x):=u_1(x)+u_2(x)$ satisfies that 
 \[ -\nabla[a_2(x) \nabla v]=\nabla[(a_1(x)-a_2(x))\nabla u_1]+b(x,u_1)-d_1u_1-(d_2+\gamma)u_2+\sigma u_3.\]
Then we can get by the energy estimate that 
 \[ \int_{\Omega}|\nabla v|^2dx+\int_{\Omega} v^2dx \leq C(\varepsilon)+2\varepsilon\int_{\Omega}u_3^2dx.\]
 From Eq.(4.3) we have by using Cauchy-Schwarz's inequality  that
 \begin{eqnarray*}
  & & a_0\int_{\Omega}|\nabla u_3|^2dx+(d_3+\sigma)\int_{\Omega}u_3^2dx\\
  & & \leq \gamma \int_{\Omega}u_2u_3dx
  \leq \frac{d_3+\sigma}{2}\int_{\Omega} u_3^2dx+\frac{\gamma}{2(d_3+\sigma)}\int_{\Omega} u_2^2dx,
 \end{eqnarray*}
 which implies 
 \[ \int_{\Omega}|\nabla u_3|^2dx+\int_{\Omega}u_3^2dx\leq C\int_{\Omega}u_2^2dx.\]
 %where $C$ depends only on known data.
 
 Now we can combine the above estimates for $u_1, v$ and $u_3$ and choose $\varepsilon$ sufficiently small to conclude
 \begin{eqnarray}
  \sum_{i=1}^{4}||\nabla u_i||_{L^{2}(\Omega)}+\sum_{i=1}^{4}\int_{\Omega}u_i^2dx\leq C. \scl
 \end{eqnarray}
% where $C$ depends only on known data.
 
 To derive a further a priori estimate, we use the Campanato estimate for elliptic equations (\cite{T1987} ) to obtain that $u_i\in C^{\alpha}(\bar{\Omega})$ and
\[ \sum_{i=1}^{4}||u_i||_{ C^{\alpha}(\bar{\Omega})}\leq C.\]
%where $C$ depends only on known data.
 
  With the above a priori estimates, we can use the Schauder's fixed-point theorem (\cite{GT1987}) to obtain the existence of a weak solution for the system (4.1)-(4.5) and the weak solution is in the space $W^{1,2}(\Omega)\bigcap C^{\alpha}(\bar{\Omega})$. We skip this step here.
   \hfill Q.E.D.
 \ \\
 {\bf Remark 4.1} The uniqueness is not expected in general since one can see that there are many nontrivial constant solutions when
 $b_1(x,s), g_1, g_2$ have the special forms as stated in the introduction.

\medskip
Next, we shall consider the steady-state solutions to the system (4.1)-(4.5).
Let 
$Z^*(x)=(u_1^*(x),u_2^*(x),u_3^*(x),u_4^*(x))$ 
be such a steady-state solution. 
For $\z>0$, we consider a small perturbation near $Z^*(x)$ and set
\[ Z(x,t)=Z^*(x)+\z Z_1(x,t), \h (x,t)\in Q,\]
where
$Z_1(x,t)=(U_1(x,t), U_2(x,t), U_3(x,t), U_4(x,t))$, with 
$U_i(x,t)=u_i(x,t)-u_i^*(x)$ for $i=1,2,3,4$.

A direct calculation shows that $Z_1$ satisfies the following linear system:
\setcounter{section}{4}
\setcounter{local}{6}
\begin{eqnarray}
 & & U_{1t}-\nabla\cdot [a_1\nabla U_1]=F_1(Z_1), \h (x,t)\in Q,\\
 & & U_{2t}-\nabla\cdot [a_2\nabla U_2]=F_2(Z_1), \h (x,t)\in Q,\scl \\
& &  U_{3t}-\nabla\cdot [a_3\nabla U_3]=F_3(Z_1), \h (x,t)\in Q,\scl\\
 & & U_{4t}-\nabla\cdot [a_4\nabla U_4]=F_4(Z_1), \h (x,t)\in Q,\scl 
 \end{eqnarray}
 subject to the initial and boundary conditions:
 \begin{eqnarray}
 & & Z_1(x,0)=Z_1(x,0), \h x\in \Omega, \scl \\
 & & \nabla_{\nu}Z_1(x,t)=0,           \h (x,t)\in \partial \Omega\times (0,\infty),\scl
 \end{eqnarray}
where the right-hand sides of the system (4.6)-(4.9) are given by 
\begin{eqnarray*}
F_1(Z_1) & = & [b_s(x,u_1^*)-\beta_1u_2^*-\beta_2h_1(u_4^*)-d_1]U_1-\beta_1u_1^*U_2+\sigma U_3-(\beta_2u_1^*h_1^{'}(u_4^*)U_4,\\
F_2(Z_1) & = & (\beta_1u_2^*+\beta_2 h_1(u_2^*0)U_1+[\beta_1u_1^*-(d_2+\gamma)]U_2+\beta_2u_1^*h_1^{'}(u_4^*)U_4,\\
F_3(Z_1) & = & \gamma U_2-(d_3+\sigma)U_3,\\
F_4(Z_1) & = & \xi U_2-h_{2s}(x,u_4^*)U_4.
\end{eqnarray*}

\noindent 
{\bf Theorem 4.1} Under the assumptions H(4.1), the steady-state solution $Z^*(x)$ to the system (4.1)-(4.5) is asymptotically stable if the following conditions hold:
\[ d_1-B_0>0, ~~d_4-G_0>0,\]
and  $\beta_1$ is suitably small,
where $B_0$ and $G_0$ are given by 
\[  B_0=\max_{x\in \Omega}|b_s(x,u_1^*)|, ~~G_0=\max_{\Omega}|h_{2s}(x,u_4^*)|.\]
{\bf Proof}. For any positive integer $k$, 
we multiply Eq.(4.1) by $U_1^k$ and integrate over $\Omega$ to obtain
\begin{eqnarray*}
& & \frac{1}{k+1}\frac{d}{dt}\int_{\Omega}U_1^{k+1}dx+\frac{4a_0}{(k+1)^2}\int_{\Omega}
|\nabla U_1^{\frac{k+1}{2}}|^2dx\\
& & +\int_{\Omega}[d_1+\beta_1u_2^*+\beta_2 h_1(u_4^*)-b_s(x,u_1^*)]U_1^{k+1}dx\\
& & \leq |J|,
\end{eqnarray*}
where $J$ is given by 
\[ J=-\beta_1\int_{\Omega}u_1^*U_2U_1^k dx +\sigma\int_{\Omega}U_3 U_1^k dx -\beta_2\int_{\Omega}u_1^*h_1^{'}(u_4^{*}) U_4U_1^k dx:=J_1+J_2+J_3.\]
Let $U_0=\max_{\Omega}u_1^*(x)$, then we can use the Young's inequality to readily get 
\begin{eqnarray*}
|J_1|&\leq& \beta_1U_0\int_{\Omega}\left[\frac{ k}{k+1}U_1^{k+1}+\frac{1}{(k+1)} U_2^{k+1}\right] dx,\\
|J_2| & \leq & \sigma\int_{\Omega}\left[\frac{ k}{k+1}U_1^{k+1}+\frac{1}{(k+1)} U_3^{k+1}\right] dx,\\
|J_3| & \leq & \beta_2U_0G_0\int_{\Omega}\left[\frac{ k}{k+1}U_1^{k+1}+\frac{1}{(k+1)} U_4^{k+1}\right] dx.
\end{eqnarray*}

Now we can easily see for sufficiently small 
$\sigma, \beta_1, \beta_2 $ that  
\begin{eqnarray*}
& & \frac{1}{k+1}\frac{d}{dt}\int_{\Omega}U_1^{k+1}dx+\frac{4a_0}{(k+1)^2}\int_{\Omega}
|\nabla U_1^{\frac{k+1}{2}}|^2dx\\
& & +[d_1+\beta_1u_2^*+\beta_2 h_1(u_4^*)-b_s(x,u_1^*)]\int_{\Omega}U_1^{k+1}dx\\
& & \leq \frac{C}{(k+1)}\int_{\Omega}\left[U_2^{k+1}+U_3^{k+1}+U_4^{k+1}\right]dx.
\end{eqnarray*}
%where $C$ is independent of $k$.

We can apply the same argument above for $U_2,U_3, U_4$ from Eq.(4.2), Eq.(4.3) and Eq.(4.4), respectively, to obtain
\begin{eqnarray*}
& & \frac{1}{k+1}\frac{d}{dt}\int_{\Omega}U_2^{k+1}dx+\frac{4a_0}{(k+1)^2}\int_{\Omega}
|\nabla U_2^{\frac{k+1}{2}}|^2dx+(d_2+\gamma-\beta_1U_0)\int_{\Omega}U_2^{k+1} dx\\
& & \leq \frac{C}{(k+1)}\int_{\Omega}\left[U_1^{k+1}+U_4^{k+1}\right]dx;\\
& & \frac{1}{k+1}\frac{d}{dt}\int_{\Omega}U_3^{k+1}dx+\frac{4a_0}{(k+1)^2}\int_{\Omega}
|\nabla U_3^{\frac{k+1}{2}}|^2dx+(d_3+\sigma)\int_{\Omega}U_3^{k+1} dx\\
& & \leq \frac{C}{(k+1)}\int_{\Omega}U_2^{k+1} dx;\\
& & \frac{1}{k+1}\frac{d}{dt}\int_{\Omega}U_4^{k+1}dx+\frac{4a_0}{(k+1)^2}\int_{\Omega}
|\nabla U_4^{\frac{k+1}{2}}|^2dx+(d_4-G_0)\int_{\Omega}U_4^{k+1} dx\\
& & \leq \frac{C}{(k+1)}\int_{\Omega}U_2^{k+1}dx.
\end{eqnarray*}
%where $C$ is independent $k$.

We now look at the quantity 
\[ Y(t)=\int_{\Omega}\left[ U_1^{k+1}+U_2^{k+1}+U_3^{k+1}+U_4^{k+1}\right]dx.\]
Noting from the assumption H(4.1) that there exists a small number, denoted by $\beta_0$, such that
\[d_1-B_0\geq \beta_0, ~~d_2+\gamma-\beta_1U_0\geq \beta_0, ~~d_3+\sigma>\beta_0, 
~~d_4-G_0\geq \beta_0, \]
we can add up the above estimates for $U_i^{k+1}$ to derive for sufficiently large $k$ that 
\[ \frac{1}{k+1}Y'(t)+\beta_0Y(t)\leq 0.\]
%where $\beta_0>0$ depends only on the known parameters.
This readily implies 
\[ Y(t)\leq C (k+1)Y(0).\]
%where $C$ is independent of $k$.
Taking the $k^{th}$-root on both sides, we obtain as 
$k\rightarrow \infty$ that
\[\sum_{i=1}^{4}\sup_{0<t<\infty}|U_i|_{L^{\infty}(\Omega)}\leq \sup_{\Omega}|Z_1(x,0)|_{L^{\infty}(\Omega)}.\]
This implies that the solution $Z_1(x,t)$ is asymptotically stable near the steady-state solution $Z^*(x)$.
\hfill Q.E.D. 

\section{Further Stability Analysis}

\setcounter{section}{5}
\setcounter{local}{1}
In this section we investigate the stability of constant steady-state solutions corresponding to the system (1.1)-(1.4). To illustrate the method and physical meaning, we further assume that  the diffusion coefficients and the death rate are constants:

\ \\
{\bf H(5.1)}. (a)  Let $a_i$ and $d_i$ be positive constants, and 
\[ a_0=min\{a_1, a_2, a_3, a_4\}, ~~d_0=min \{d_1, d_2, d_3, d_4\}.\]

(b) Functions $b$, $h_1$ and $h_2$ are of the following 
forms for two constants $b_0$ and $g_0$: 
\[ b(x,t,s)=b_0s(1-\frac{s}{k_{1}}), ~~h_1(s)=\frac{s}{s+k_2}, ~~h_2(s)=g_0s(1-\frac{s}{k_3}).\]

Consider the corresponding steady-state system in $\Omega$:
\setcounter{section}{5}
\setcounter{local}{1}
\begin{eqnarray}
 -\nabla\cdot [a_1(x)\nabla u_1] & = & b(x,u_1)-\beta_1u_1u_2-\beta_2u_1\cdot h_1(u_4)-d_1u_1+\sigma u_3,\\
 -\nabla\cdot [a_2(x)\nabla u_2] & = & \beta_1u_1u_2+\beta_2u_1 \cdot h_1(u_4)-(d_2+\gamma)u_2,\scl \\
 -\nabla\cdot [a_3(x)\nabla u_3] & = & \gamma u_2-(d_3+\sigma)u_3,\scl \\
 -\nabla\cdot [a_4(x)\nabla u_4] & = & \xi u_2 +h_2(x,u_4)-d_4 u_4\scl
 \end{eqnarray}
 subject to the boundary condition
 \begin{eqnarray}
 \partial_{\nu}U(x) =0, \h x\in \partial \Omega,\scl
 \end{eqnarray}
 where  $U(x)=(u_1(x), u_2(x), u_3(x), u_4(x))$.

We can easily derive from (5.1) to (5.4) that 
\[ (d_1-b_0)\int_{\Omega}u_1dx+d_2\int_{\Omega}u_2dx+d_3\int_{\Omega}u_3dx+\frac{b_0}{K_1}\int_{\Omega}u_1^2dx=0.\]
\[(d_4-g_0)\int_{\Omega}u_4dx+\frac{g_0}{K_2}\int_{\Omega}u_4^2dx=\xi\int_{\Omega}u_2dx, \]
from which we readily see that there exists one trivial solution, i.e., $u_1=u_2=u_3=u_4=0$ if $b_0\leq d_1$ and $g_0\leq d_4$.

On the other hands, we can also see that there are two sets of steady-state solutions. The first set of constant solutions requires
$b_0> d_1$ and $g_0> d_4$:
\begin{eqnarray*}
& & Z_1=(0,0,0,0); ~~Z_2=(\frac{K_1(b_0-d_1)}{b_{0}}, 0,0, 0); 
~~Z_3=(0,0,0,\frac{K_2(g_0-d_4)}{g_{0}}).
\end{eqnarray*}
There exists another set of constant solutions:
\[ Z_4=\left\{(S,I,R,B): R=\frac{\gamma}{d_3+\sigma}I. \right\}, \]
where $S,I$ and $B$ are the solutions of the following nonlinear system:
\begin{eqnarray}
& & \frac{b_0}{K_1}S^2-(b_0-d_1)S+\left(d_2+\gamma-\frac{\sigma \gamma}{d_3+\sigma}\right)I=0,\\
& & \frac{g_0}{K_2}B^2-(g_0-d_4)B-\xi I=0,\scl\\
& & S=\frac{(d_2+\gamma)I}{\beta_1I+\beta_2 h_1(B)}.\scl
\end{eqnarray}

\ \\
{\bf Lemma 5.1}. The nonlinear system (5.5)-(5.7) has at least one solution if and only if the following condition holds:
\[\frac{K_1(b_0-d_1)}{2b_0}>\frac{d_2+\gamma}{\beta_2}.\]
\ \\
{\bf Proof}: We first derive a necessary condition which will ensure the existence of 
a nontrivial constant solution.
By solving the quadratic equation (5.5) for $S$, we obtain
\begin{eqnarray*}
& & S_1=\frac{(b_0-d_1)+ \sqrt{(b_0-d_1)^2-\frac{4b_0}{K_1}[(d_2+\gamma)-\frac{\sigma \gamma}{d_3+\sigma}]I}}{\frac{2b_0}{K_1}},\\
& & S_2=\frac{(b_0-d_1)- \sqrt{(b_0-d_1)^2-\frac{4b_0}{K_1}[(d_2+\gamma)-\frac{\sigma \gamma}{d_3+\sigma}]I}}{\frac{2b_0}{K_1}}.
\end{eqnarray*}
Noting that 
\[ (d_2+\gamma)-\frac{\sigma\gamma}{d_3+\sigma}>0, \]
we see that the range of $I$ must satisfy
\[ 0\leq I\leq I^*:=\frac{K_1(b_0-d_1)^2}{4b_0[(d_2+\gamma)-\frac{\sigma\gamma}{d_3+\sigma}]}.\]
But we can see from Eq.(5.7) that 
\[S=\frac{(d_2+\gamma)I}{\beta_1I+\beta_2 h_1(B)}=\frac{d_2+\gamma}{\beta_1}[1-\frac{\beta_2 h_1(B)}{\beta_1I+\beta_2h_1(B)}].\]
If we consider $S$ as a function of $I$, i.e., $S=S(I)$, we get 
\[ S(0)=0, ~~S'(I)>0, ~~S(\infty)=\frac{d_2+\gamma}{\beta_1}.\]
On the other hand, if we consider $S_1$ as a function of $I$, i.e.,  $S_1=S_1(I)$, 
then we have 
\[ S_1(0)=\frac{K_1(b_0-d_1)}{b_0}, ~~S_1'(I)<0.\]
We readily see that
\[ \min_{I\in[0,I^*]}S_1(I)=S_1(I^*)=\frac{K_1(b_0-d_1)}{2b_0}, ~~\max_{I\in [0,I^*]}S_1(I)=S_1(0)=\frac{K_1(b_0-d_1)}{b_0}.\]
Consequently, $S(I)$ and $S_1(I)$  have an intersection point if and only if
\[ \frac{K_1(b_0-d_1)}{2b_0}>\frac{d_2+\gamma}{\beta_2}.\]
Moreover, the intersection point is unique since both $S(I)$ and $S_1(I)$ are monotone functions.

Similarly, we see for $S_2$, 
\[ S_2(0)=0, ~~S_2'(I)>0, ~~S_2''(I)>0.\]
Hence we have 
\[\max_{I\in [0,I^*]}S_2(I)=\frac{K_1(b_0-d_1)}{2b_0}\]
The above indicates the existence of an intersection point between $S(I)$ and $S_2(I)$ as long as
\[\frac{K_1(b_0-d_1)}{2b_0}>\frac{d_2+\gamma}{\beta_2}.\]
Once $I$ and $S$ are determined, one can easily solve for $B$ from Eq.(5.6):
\[ B=\frac{K_2\left[(g_0-d_4)+\sqrt{(g_0-d_4)^2+\frac{4g_0\xi}{K_2}I}\,\right]}{2g_0}.\]

\hfill Q.E.D.

\ \\
{\bf Proof of Theorem 2.2}.  
Let $A$ be the diagonal matrix with the diffusion coefficients $a_i$.
We can calculate the Jacobian matrix for the nonlinear reaction terms from system (2.1)-(2.4):
\[B_1(Z)= \left( \frac{\partial f_i}{\partial u_{i}} \right)_{4\times 4}.\] 

For $Z_1=(0,0,0,0)$, it is easy to see the $4 \times 4$~matrix:
\[B_1(Z_1)= \left( \begin{array}{cccc}
b_0-d_1 & 0 & \sigma & -\frac{\beta_2K_1(b_0-d_1)}{b_0K_2}\\
0  &\frac{\beta_2K_1(b_0-d_1)}{b_0} -(d_2+\gamma) & 0 & \frac{\beta_2K_1(b_0-d_1)}{K_2b_0} \\
0 & \gamma & -(d_3+\sigma) & 0\\
0 & \xi & 0 & g_0-d_4
 \end{array} \right).\] 
Let  $0\leq \lambda_1<\lambda_2<\cdots $ be the eighenvalue of the Laplacian operator subject to the homogeneous Neumann boundary condition.

 It is easy to calculate the eigenvalues of $A_j(Z_1)=DF(Z_1)-\lambda_jA$:
 \[ \mu_{1j}=b_0-d_1-\lambda_j a_1, \mu_{2j}=-(d_2+\gamma)-\lambda_j a_2, \mu_{3j}=-(d_3+\sigma)-\lambda_ja_3, \mu_{4j}=
 g_0-d_4-\lambda_j a_4.\]
  Since $\lambda_1=0$ is the first eigenvalue and $b_0\geq d_1$ and $g_0\geq d_4$, it follows that
  $Z_1=(0,0,0,0)$ is unstable unless $b_0\leq d_1, g_0 \leq d_4$. 
  
  Since $\lambda_j\geq 0$, the eigenvalues indicate that the stability of $Z_1$ is not affected by the diffusion processes.
This is clear since the birth rate is greater than the death rate. The population must be positive for a long time.
  
  For $Z_2=(\frac{K_1(b_0-d_1)}{b_{0}}, 0,0, 0)$, we can see the $4 \times 4$~matrix:
  \[B_1(Z_2)= \left( \begin{array}{cccc}
-(b_0-d_1) & -\frac{K_1\beta_1(b_0-d_1)}{b_0} & \sigma & -\frac{\beta_2K_1(b_0-d_1)}{b_0K_2} \\
0  & \frac{\beta_1K_1(b_0-d_1)}{b_0}-(d_2+\gamma) & 0 & \frac{\beta_2K_1(b_0-d_1)}{b_0K_2} \\
0 & \gamma & -(d_3+\sigma) & 0\\
0 & \xi & 0 & g_0-d_4
 \end{array} \right).\] 
Then we consider   
\[ A_j(Z_2)=DF(Z_2)-\lambda_jA,\]
and see its characteristic polynomial, denoted by $P(\mu)$, is equal to
\begin{eqnarray*}
P(\mu)= & & (b_0-d_1-\lambda_ja_1-\mu)(d_3+\sigma+\lambda_ja_3+\mu)\\
& &\{ [\mu^2-[(g_0-d_4-\lambda_ja_4+m_0-(d_2+\gamma+\lambda_j a_2)\mu\\
& & + [m_0-(d_2+\gamma+\lambda_j a_2)][g_0-d_4-\lambda_j a_4]-\xi m_0\}.
\end{eqnarray*}
where
\[ m_0=\frac{\beta_2K_1(b_0-d_1)}{b_0}.\]
We obtain the eigenvalues
\begin{eqnarray*}
\mu_1 & = & -(b_0-d_1)-\lambda_ja_1,\\
\mu_2 & = & -(d_3+\sigma+\lambda_ja_3),\\
\mu_3 & = & \frac{M_1+\sqrt{M_1^2-4M_2}}{2},\\
\mu_4 & = & \frac{M_1-\sqrt{M_1^2-4M_2}}{2},
\end{eqnarray*}
where 
\begin{eqnarray*}
M_1 & = & m_0-(d_2+\gamma+\lambda_ja_2)+(g_0-d_4-\lambda_ja_4);\\
M_2 & = & [m_0-(d_2+\gamma+\lambda_ja_2)][g_0-(d_4+\lambda_ja_4)]-\xi m_0
\end{eqnarray*}
It follows that
$Z_2$ is locally stable if $M_1<0$ and $M_2>0$ and $Z_2$ is unstable for either $M_1>0$ or $M_2<0$
or $M_1^2-4M_2>0$ when $M_2>0$.
On the other hand, we know 
\[ \lambda_j\rightarrow \infty ~ \mbox{as $j\rightarrow \infty$},\]
and $M_1^2-4M_2>0$.
Consequently, we conclude that $Z_2$ is an unstable steady-state solution.

Now we calculate $A_j(Z_3)$:
\[ A_j(Z_3)=DF(Z_3)-\lambda_j A.\]
For $Z_3=(0,0,0,\frac{K_2(g_0-d_4)}{g_{0}})$, we can see the $4 \times 4$~matrix:
  \[B_1(Z_3)= \left( \begin{array}{cccc}
(b_0-d_1) & 0 & \sigma & 0 \\
\frac{\beta_2(g_0-d_4) }{(2g_0-d_4)}  & -(d_2+\gamma) & 0 & 0 \\
0 & \gamma & -(d_3+\sigma) & 0\\
0 & \xi & 0 & -( g_0-d_4)
 \end{array} \right).\] 
We know the characteristic polynomial for the matrix $A_j(Z_3)=DF(Z_3)-\mu I_{4\times 4}$ is equal to
\begin{eqnarray*}
P(\mu)= & & |A_j(Z_3)|=-[(g_0-d_4+\lambda_j a_4)+\mu]P_0(\mu),
\end{eqnarray*}
where
\[ P_0(\mu)=\ [(b_0-d_1-\lambda_j a_1-\mu)(d_3+\sigma+\lambda_ja_3+\mu)(d_2+\gamma+\lambda_ja_2+\mu)+\frac{\sigma \gamma \beta_2(g_0-d_4)}{2g_0-d_4}.\]

Hence, the first eigenvalue is equal to
\begin{eqnarray*}
\mu_1 & = & -(g_0-d_1+\lambda_j a_4),
\end{eqnarray*}
To see the rest of eigenvalues of $P(\mu)$, we use a lemma from Yin-Chen-Wang \cite{YCW2017}.
\ \\
{\bf Lemma 5.2} Let $p>0$, $q$ and $h$ be constants, and 
\[ P_0(\mu)=\mu^3+p\mu^2+q\mu +h=0.\]
Then it holds that \\
(a) If $h<0$, there exists a positive root;\\
(b) If $0<h<pq$, all roots have negative real parts;\\
(c) If $pq<h$, there is a root with positive real part;\\
(d) If $pq=h$, the roots are $\mu_1=-p, \mu_2=\sqrt{-q}, \mu_3=-\sqrt{-q}.$

Let 
\[ P_0(\mu)=\mu^3+p\mu^2+q\mu+h,\]
with its coefficients given by 
\begin{eqnarray*}
p & = & (d_2+\gamma+\lambda_j a_2)+(d_3+\sigma+\lambda_ja_3)-(b_0-d_1-\lambda_ja_1);\\
q & = & (d_2+\gamma+\lambda_j a_2)(d_3+\sigma+\lambda_j a_3)-(b_0-d_1-\lambda_j a_1)[(d_2+\gamma+\lambda_j a_2)+(d_3+\sigma+\lambda_ja_3)];\\
h & = & (d_1+\lambda_ja_1-b_0)(d_2+\gamma+\lambda_ja_2)(d_3+\sigma+\lambda_j a_3)-\frac{\sigma\gamma\beta_2(g_0-d_4)}{2g_0-d_4}.
\end{eqnarray*}

Since $\lambda_1=0$ is one of the eigenvalues and $d_1-b_0<0, g_0-d_4>0$, we see $h<0$ from the expression of $h$, so $Z_3$ is unstable.

Finally, we study the stability of $Z_4$. Since $u_4$ always has positive solutions as long as $u_2$ is positive, it does not affect
the stability of other variables. We only need to focus on the stability of $(u_1,u_2,u_3)$. Furthermore, since $\lambda_1=0$ is the first eigenvalue, the rest of eigenvalues have the same sign with $d_i$ which increases the stability of the solution. Therefore, we only need to find the conditions for the stability  when $\lambda_1=0$.

It is easy to calculate the Jacobian matrix 
\[B_1^*= \left( \begin{array}{ccc}
-L_0 & -\beta_1S_0 & \sigma \\
\beta_1I_0+\beta_2h(B_0) & -(d_2+\gamma) & 0 \\
0 & \gamma & -(d_3+\sigma) 
 \end{array} \right)\] 
 where
 \[ L_0= (d_1-b_0)+\frac{2b_0S_0}{K_1}+\beta_1I_0+\beta_2h_1(B_0).\]
The characteristic polynomial of $B_1^*$ is equal to
\[ P(\mu)=\mu^3+p_0\mu^2+q_0\mu+h_0=0.\]
where
\begin{eqnarray*}
p_0 & = & L_0+(d_2+\gamma)+(d_3+\sigma)+L_0;\\
q_0 & = & (d_3+\sigma)(L_0+d_2+\gamma)+L_0(d_2+\gamma)
+\beta_1S_0(\beta_1I_0+\beta_2h_1(B_0));\\
h_0 & = & (d_3+\sigma)[L_0(d_2+\gamma)+\beta_1S_0(\beta_1I_0+\beta_2h_1(B_0)]
-\sigma\gamma(\beta_1 I_0+\beta_2h_1(B_0)).
\end{eqnarray*}
By Lemma 5.2, we can see the stability or instability of the steady-state solution precisely when parameters varies. 
In particular, when $L_0>0$, if $\sigma, \gamma, \beta_1$ and $\beta_2$ are sufficiently small, we see
the condition $0<h_0<p_0q_0$ holds. Consequently, the steady-state solution $(S_0,I_0,R_0)$ is stable. This result confirms the result of Theorem 2.2 about the stability analysis of the steady-state solution. \hfill Q.E.D.

 \section{Conclusion}
 
 In this paper we have studied a nonlinear mathematical model for an epidemic caused by cholera without life-time immunity.
 The diffusion coefficients are different for each species. Moreover, these coefficients are allowed to be dependent upon the concentration as well as the space location and time. The resulting model system is strongly coupled. We established the global well-posedness for the coupled reaction-diffusion system under some very mild conditions on the given data. Moreover,
 we have analyzed the linear stability for the steady-state solutions and proved that there is a turing phenomenon when the diffusion coefficients are different. This result indicates that there are some fundamental differences between the ODE model and the corresponding PDE model.
 These results show that the mathematical model is well-defined and can be used by other researchers to conduct the field study. 
 The theoretical results obtained in this paper lays a solid foundation for other scientists in related fields to further study
 more constructive qualitative properties of the solutions. The study will provide scientists a deeper understanding of the dynamics of the interaction between bacteria and susceptible, infected and recovered species. We have used many ideas and techniques from the elliptic and parabolic equations, particularly, the energy method and Sobolev's inequalities. 
    There are some open questions that remain to be answered, and further studies are needed.
 
 \ \\
 {\bf Acknowledgements}.
   This work was motivated by some open questions raised by Professor K. Yamazaki  from Texas Tech University and Professor Jin Wang from University of Tennessee at Chattanooga in WSU biological seminar series. The authors would like to thank them for some helpful discussions about the model. 
   The work of the second author was substantially supported by Hong Kong RGC General Research Fund (projects 14306921 and 14306719).

\end{document}